\documentclass[12pt,twoside]{paper}
\usepackage{amsmath}
\usepackage{hyperref}
\usepackage{amsfonts}
\usepackage{amssymb}

\usepackage{graphics}
\usepackage[pdftex]{graphicx}

\def\bee{\begin{equation}}
\def\ee{\end{equation}}
\def\Li{{\rm Li}}
\input epsf

\begin{document}

\thispagestyle{empty}
\bigskip\bigskip
\centerline{    }
\vskip 4 cm
\centerline{\Large\bf An analog of the Skewes number for twin primes}
\bigskip\bigskip\bigskip
\centerline{\large\sl Marek Wolf}
\bigskip
\centerline{\sf Institute of Theoretical Physics, University of Wroc{\l}aw}
\centerline{\sf Pl.Maxa Borna 9, PL-50-204 Wroc{\l}aw, Poland}
\centerline{\small e-mail:  mwolf@ift.uni.wroc.pl}

\bigskip\bigskip\bigskip

\begin{center}
{\bf Abstract}\\
\end{center}

\begin{minipage}{12.8cm}
The results of the computer investigation of the sign changes of the
difference between the number of twin primes $\pi_2(x)$ and the
Hardy--Littlewood  conjecture  $c_2\Li_2(x)$ are reported. It turns out that
$\pi_2(x) - c_2\Li_2(x)$  changes the sign at unexpectedly low values of $x$ and
for $x<2^{42}$ there are over 90000 sign changes of this difference.
It is conjectured that the number of sign changes of
$\pi_2(x) - c_2\Li_2(x)$ for $x\in (1, T)$   is given by $\sqrt T/\log(T)$.
\end{minipage}

\bigskip\bigskip\bigskip

\vfill
\eject

\pagestyle{myheadings}

Let $\pi(x)$ be the number of primes less than $x$ and let $\Li(x)$
denote the logarithmic integral:
\begin{equation}
\Li(x)=\int_2^x \frac {du}{ \log(u)}.
\label{Li}
\end{equation}
The Prime Number Theorem tells us that $\Li(x)/\pi(x)$ tends to 1 for
$x\rightarrow\infty$
and the available data (see \cite{Odlyzko}) shows that always
$\Li(x)>\pi(x)$. This last experimental observation was the reason of
the common belief in the past, that the inequality
$\Li(x)>\pi(x)$ is generally valid.
However, in 1914 J.E. Littlewood has shown \cite{Littlewood} (see also \cite{Ellison})
that the difference between the number of primes smaller than $x$ and the
logarithmic integral up to $x$ infinitely often changes the sign.
In 1933 S. Skewes \cite{SkewesI}
assuming the truth of the Riemann hypothesis has estimated that for sure
$d(x)=\pi(x)-\Li(x)$ changes sign for some
$x_0 < 10^{10^{10^{34}}}$. The smallest value $x_0$ such that
for the first time $\pi(x_0)>\Li(x_0)$ holds is called Skewes number. In 1955
Skewes \cite{SkewesII} has found, without assuming
the Riemann hypotheses, that $d(x)$ changes sign at some
$$
x_0 < \exp\exp\exp\exp(7.705)<10^{10^{10^{10^3}}}.
$$
This enormous bound for
$x_0$ was reduced by Cohen and Mayhew \cite{Cohen} to $x_0<
10^{10^{529.7}}$ without using the Riemann hypothesis. In 1966
Lehman \cite{Lehman} has shown that between $1.53\times 10^{1165}$ and
$1.65\times 10^{1165}$ there are more than $10^{500}$ successive integers
$x$ for which $\pi(x)> \Li(x)$. Following the method of Lehman in 1987
H.J.J. te Riele  \cite{Riele} has shown that between $6.62\times 10^{370}$ and
$6.69\times 10^{370}$ there are more than $10^{180}$ successive integers
$x$ for which $d(x)>0$. The lowest present day known value of the Skewes number is around
$10^{316}$, see \cite{Bays} and \cite{Demichel}.

The number of sign changes of the difference $d(x) =\pi(x)-\Li(x)$ for
$x$ in a given interval $(1, T)$, which is commonly denoted by $\nu(T)$,
see \cite{Ellison}, was treated for the first time  by A.E. Ingham in 1935 \cite{Ingham1}
chapter V, \cite{Ingham2} and next by
S. Knapowski \cite{Knapowski}.
Regarding the number of sign changes of $d(x)$ in the interval
$(1,T)$,  Knapowski \cite{Knapowski} proved
\begin{equation}
\nu(T)\geq e^{-35}\log \log \log \log T
\label{nu}
\end{equation}
provided $T\geq \exp\exp\exp\exp(35)$. Further results about $\nu(T)$ were obtained
by  J. Pintz \cite{Pintz} and J. Kaczorowski \cite{Kaczorowski}

In this paper I am going to look for the analog of the Skewes number for the twin primes.

Let $\pi_2(x)$ denote the number of twin primes smaller than $x$. Then
the unproved (see however \cite{Rubinstein})
conjecture B of Hardy and Littlewood \cite{Hardy and Littlewood} on
the number of prime pairs $p, p+d$ applied to the case $d=2$ gives, that
\begin{equation}
\pi_2(x) \sim C_2\Li_2(x) \equiv C_2 \int_a^x \frac{u}{\log^2(u)} du,
\label{conj}
\end{equation}
where $C_2$ is called ``twin constant'' and is defined by the following
infinite product:
\begin{equation}
C_2 \equiv 2 \prod_{p > 2} \biggl( 1 - \frac{1}{(p - 1)^2}\biggr) =
1.32032\ldots
\label{stalac2}
\end{equation}

Usually the lower limit of integration $a$ in (\ref{conj}) is chosen 2,
but the author believes that the proper choice for the lower limit of
integration should be 5, not 2, because (3,5) is the first twin pair.
(Analogously in (\ref{Li}) the lower limit of integration is 2, to
ensure that $\Li(2)=0$.)

For the first time the conjecture (\ref{conj}) was checked computationaly by R. P. Brent
\cite{Brent}. This author noticed the sign changes of the difference
$\pi_2(x)-C_2 {\rm Li_2}(x)$ but did not elaborate about this further.
I have looked on the difference $d_2(x) = \pi_2(x)-C_2\Li_2(x)$ using the
computer  for $T$ up to $2^{42}\approx 4.4\times 10^{12}$. Like for usual primes initially
$C_2 {\rm Li}_2(x)> \pi_2(x)$, but
surprisingly, it turns out that there is
a lot sign changes of $d_2(x)=\pi_2(x)-C_2\Li_2(x)$ for $x$ in the range
$(1, 2^{42})$. For the case of the lower integration limit $a=2$ in (\ref{conj}) the first
sign change of
$d_2(x)$ appears at the twin pair (1369391, 1369393). However, for the
choice $a=5$ the first sign change of $d_2(x)$ appears already at the
pair (41, 43)! Next ``Skewes'' pairs  for $a=5$ are (6959, 6961) and (7127, 7129). After
that there is a gap in the sign changes up to (1353257, 1353259), what
is comparable to the case of $a=2$. Such a behavior is reasonably, because the
difference between $a=2$ and $a=5$ is relevant only at low $x$
and the contribution stemming from the interval (2,5) is becoming
negligible for large $x$.

Let $\nu_2(T)$ denote, by analogy with usual primes,
the number of sign changes of $d_2(x)$ in the interval
$(1,T)$.   The Table I contains the recorded number of
sign changes of $\pi_2(x) - C_2\Li_2(x)$ up to $T=2^{22}, 2^{23}, \ldots,
2^{42}$ for both choices of the parameter $a$.
The values of $T$ searched by the direct checking are of small magnitude from
the point of view of mathematics, but large for modern computers.
The observed numbers $\nu_2(T)$ behave
very erratically, see Fig.1, in particular there are large gaps without
any change of sign of the $d_2(x)$.  If one assumes the power-like dependence of $\nu_2(T)$
then the fit by the least square method  gives  the function  $\alpha T^\beta$,
where $\alpha=0.0766741$ and $\beta=0.479031\ldots$. Instead of such  accidentally
looking parameters  for the pure power-like behavior  (especially $\alpha$ has a very
small value) after a few trials I have picked the function $\sqrt T/\log(T)$ as an
approximation to $\nu_2(T)$ as a much nicer function without any random parameters
involved and simultaneously taking values very close to the fit $\alpha T^\beta$,
see Figure 1. Thus we state the conjecture:
\begin{equation}
\nu_2(T) = \sqrt T/\log(T)~+ ~error~~term
\label{conjecture}
\end{equation}
Let us stress in favor of
(\ref{conjecture}) that there are 7 crosses of the curve $\sqrt T/\log(T)$
with the staircase-like plot of $\nu_2(T)$ obtained directly from the computer data.
The last
column  in the Table 1 contains the values of the
function $\sqrt T/\log(T)$. If the conjecture (\ref{conjecture})
is true, then there is infinity of twins. Also if (\ref{conjecture}) is
valid, it means that the estimation (\ref{conj}) is in some sense more
accurate than (\ref{Li}), because there are more points, where
(\ref{conj}) exactly reproduces the actual number of twins --- for
(\ref{Li}) there is much less such values of $x$ that $\Li(x)$ is equal
to $\pi(x)$, see (\ref{nu}). However, presumably the (unknown) error term
in (\ref{conj}) is larger than error term for $\pi(x)$.

The difference of many  hundreds of orders between  values of $x$ such that
$\pi(x)- {\rm Li}(x)$ and  $\pi_2(x)-C_2{\rm Li}(x)$ changes the sign for the  first
time seems to be very astonishing. Let me give an example from physics:
the energy of the ground states of the hydrogen and helium are respectively
-13.6 eV and  -79 eV and do not differ by hundreds of orders!

I have tested the numerical results using several different computers, programs
and compilers. In particular, to calculate the integral ${\rm Li}_2(x)$ I have used the
8--point self--adaptive Newton--Cotes method and the 10--point Gauss
method. This integral was calculated numerically in successive intervals between
consecutive twins and added to the previous value. It seems to be natural that different
methods gave exactly the same values of $\Li_2(x)$
since the integrand in (\ref{conj}) is a very well behaved function. At
least up to $T=2^{32}$ (this limitation stems from the fact that some compilers did
not allow larger integers than $2^{32}$) all results obtained by different runs were
exactly the same.

\vfill
\eject

\newpage

\vskip 0.4cm
\begin{table}
\renewcommand{\arraystretch}{1}
\begin{center}
{\sf TABLE {\bf 1}}\\
{\small The number of sign changes of $d_2(x)$.
The case $a=2$ is in second column, and third column contains data for $a=5$, while
the last column presents values obtained from (\ref{conjecture}).}\\
\bigskip
\begin{tabular}{||c|c|c||c||} \hline
$T$ & $ \nu_2(T) $ for $ a=2 $ & $ \nu_2(T)$ for $a=5$ & $\frac{\sqrt T}{\log(T)} $ \\ \hline
$2^{22}$ &           29 &            32 &       134\\ \hline
$2^{23}$ &           29 &            32 &       182\\ \hline
$2^{24}$ &           29 &            32 &       246\\ \hline
$2^{25}$ &           29 &            32 &       334\\ \hline
$2^{26}$ &          238 &           269 &       455\\ \hline
$2^{27}$ &          854 &           942 &       619\\ \hline
$2^{28}$ &         1226 &          1401 &       844\\ \hline
$2^{29}$ &         1226 &          1401 &      1153\\ \hline
$2^{30}$ &         1226 &          1401 &      1576\\ \hline
$2^{31}$ &         1226 &          1401 &      2157\\ \hline
$2^{32}$ &         2854 &          3045 &      2955\\ \hline
$2^{33}$ &         7383 &          7358 &      4052\\ \hline
$2^{34}$ &         9115 &          8974 &      5562\\ \hline
$2^{35}$ &        12682 &         12431 &      7641\\ \hline
$2^{36}$ &        23634 &         23103 &     10505\\ \hline
$2^{37}$ &        31641 &         30770 &     14455\\ \hline
$2^{38}$ &        31641 &         30770 &     19905\\ \hline
$2^{39}$ &        31641 &         30770 &     27428\\ \hline
$2^{40}$ &        38899 &         37904 &     37819\\ \hline
$2^{41}$ &        55106 &         54179 &     52180\\ \hline
$2^{42}$ &        90355 &         89768 &     72037\\ \hline
\end{tabular} \\
\end{center}
\end{table}

\newpage

\begin{figure}
\includegraphics[width=12cm,angle=0, scale=1.0]{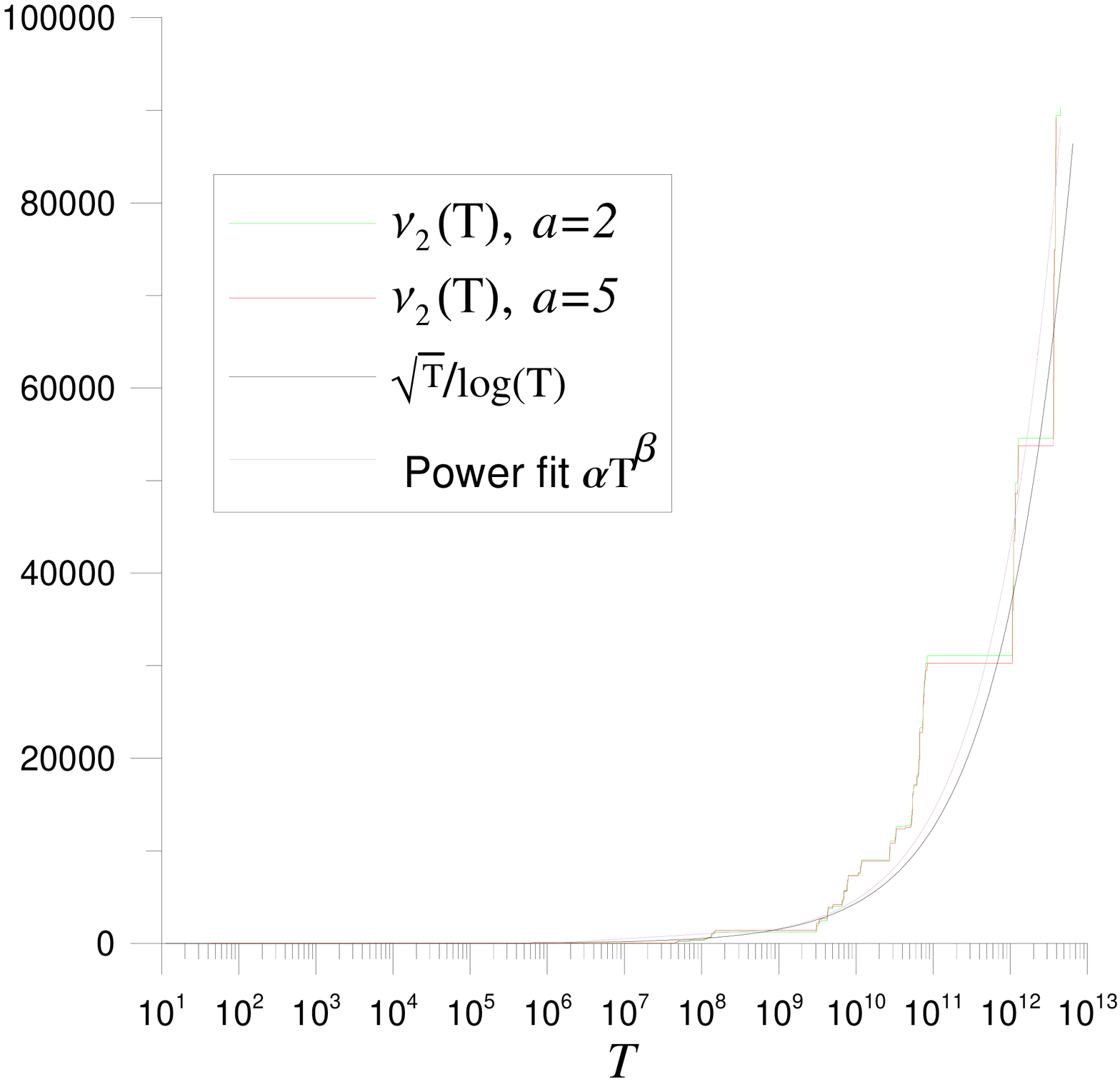} \\
Fig.1 The plot showing the comparison of the actual values of $\nu_2(T)$ for
$a=2$ and $a=5$ found by a computer search with the conjecture  (\ref{conjecture}).
\end{figure}

\end{document}